\documentclass{article}
\usepackage{amssymb}
\usepackage{amsfonts}
\begin{document}

\title{On the increase of Gronwall function value at the multiplication of its argument by a prime}
\author{Aleksandr Morkotun\thanks{E-mail: a.morkotun@gmail.com}}

\huge
\begin{center}
On the increase of Gr\"onwall function value at the multiplication of its argument by a prime
\end{center}

\par\bigskip

\Large
\begin{center}
Aleksandr Morkotun
\end{center}

\large
\begin{center}
E-mail: a.morkotun@gmail.com
\end{center}

\par\bigskip
\par\bigskip

\normalsize

\noindent {\bf Abstract.} We consider the function $G(n)=\frac{\sigma(n)}{n\log\log n}$ (where $\sigma(n)=\sum_{d|n}d$) and set an imposed condition on its argument $n$, the fulfillment of which is sufficient for the existence of a prime $p$, at which $G(np)>G(n)$. This inequality is of interest in connection with the Robin's inequality. The paper also presents the results of numerical experiment conducted with superabundant numbers.

\par\bigskip


\large

\par\bigskip
\par\bigskip

\section{Introduction}

The well-known Robin's theorem~\cite{Robin} proclaims the equivalence of the Riemann hypothesis of non-trivial zeros of the Riemann zeta function and elementary statement that the inequality
\begin{equation}
\label{1}
G(n)<e^\gamma
\end{equation}
is true at any $n>5040$. Function $G$, featured here, is called the Gr\"onwall function~\cite{Gronwall}. It is defined by the equality
$$
G(n)=\frac{\sigma(n)}{n\log\log n},
$$
in which $\sigma(n)$ denotes the sum of divisors of $n$ and $\gamma$ --- the Euler-Mascheroni constant.

Many publications are devoted to the study of function $G$ behavior. Their review can be found in~\cite{Nicolas-1} and ~\cite{NazYak}. The work~\cite{Nicolas-2} contains the following remarkable theorem.

\par\medskip
\noindent {\bf Theorem 1.} {\itshape Riemann hypothesis is true if and only if for each $n>$~$5040$ there exists such a positive integer $m$, at which $G(mn)>G(n)$.}
\par\medskip

Our main result is a theorem, which sets an imposed condition on $n$ sufficient for the existence of a {\it prime} factor $m$:

\noindent {\bf Theorem 2.} {\itshape Let $n=\prod\limits_{i=1}^k p_i^{a_i}$, where $p_i$ denotes the {\it i}-th prime number ($p_1=2$, $p_2=3$ etc.). If for certain $l$ the inequality
\begin{equation}
\label{2}
p_l^{a_l+1}<\log n
\end{equation}
is true, then the inequality
\begin{equation}
\label{3}
G(np_l)>G(n)
\end{equation}
is also true.}
\par\medskip

Material of the paper is as follows. In the second section we give a proof of Theorem~2, the third section describes a numerical experiment which deals with the so-called superabundant numbers, and the last section contains concluding remarks.

\section{Proof of the Theorem 2}

We need a proposition that establishes a necessary and sufficient condition for the inequality~(\ref{3}), and an auxiliary lemma.

\par\medskip
\noindent {\bf Proposition 3.} {\itshape Let $l$ be positive integer and $n=\prod\limits_{i=1}^k p_i^{a_i} \ge 3$ ($p_i$ is $i$-th prime). Then
$$
G(np_l)>G(n) \Longleftrightarrow (\log n)^{\left(\sum\limits_{\alpha=1}^{a_l+1}p_l^\alpha\right)^{-1}}>1+\frac{\log p_l}{\log n}.
$$}
\par\medskip

\noindent {\it Proof.} Let us consider the ratio $\frac{G(np_l)}{G(n)}$. Provided that $\sigma(n)=\prod\limits_{i=1}^k \frac{p_i^{a_i+1}-1}{p_i-1}$, we have:
$$
\frac{G(np_l)}{G(n)}=\frac{\sigma(np_l)}{\sigma(n)}\cdot\frac{\log\log n}{p_l\cdot\log\log (np_l)}=\frac{p_l^{a_l+2}-1}{p_l^{a_l+2}-p_l}\cdot\frac{\log\log n}{\log\log (np_l)}.
$$
Thus,
$$
G(np_l)>G(n)\Longleftrightarrow(\log n)^\frac{p_l^{a_l+2}-1}{p_l^{a_l+2}-p_l}>\log n+\log p_l \Longleftrightarrow
(\log n)^{\left(\sum\limits_{\alpha=1}^{a_l+1}p_l^\alpha\right)^{-1}}>1+\frac{\log p_l}{\log n}.
\hspace{0.5cm}\square
$$

\par\medskip
\noindent {\bf Lemma 4.} {\itshape Let $t$ and $s$ be positive integers, and $\xi$ --- real number greater than $t^s$. Then the following inequality holds:
$$
\xi^{\left(\sum\limits_{\alpha=1}^{s}t^\alpha \right)^{-1}}>1+\frac{\log t}{\xi}.
$$}
\par\medskip

\noindent {\it Proof.} Since $\xi>t^s$ and $s\left(\sum_{\alpha=1}^{s}t^\alpha \right)^{-1}\ge t^{-s}$, we have:
$$
\xi^{\left(\sum\limits_{\alpha=1}^{s}t^\alpha \right)^{-1}}>t^{s\left(\sum\limits_{\alpha=1}^{s}t^\alpha \right)^{-1}}\ge 1+\log t^{t^{-s}}=1+\frac{\log t}{t^s}>1+\frac{\log t}{\xi}.
\hspace{0.5cm}\square
$$

\noindent {\it Proof of the Theorem 2.} Provided that in Lemma 4 $\xi=\log n$, $t=p_l$, and $s=a_l+1$, we can write the following:
$$
\log n>p_l^{a_l+1}\Longrightarrow (\log n)^{\left(\sum\limits_{\alpha=1}^{a_l+1}p_l^\alpha\right)^{-1}}>1+\frac{\log p_l}{\log n}.
$$
Combining this with Proposition 3 follows the statement of the Theorem 2.
$\hspace{0.3cm}\square$

\section{Numerical experiment}

In connection with obtaining Theorem 2 it is natural to ask how common the numbers satisfying the condition~(\ref{2}) are, and how the fact of their presence may be helpful for proving the Robin's inequality~(\ref{1}) in the general case.

When studying behavior of the Gr\"onwall function so-called superabundant numbers (SA numbers) are of special interest. The positive integer $n$ is called superabundant~\cite{AE,Ramanujan}, if the inequality $\frac{\sigma(m)}{m}<\frac{\sigma(n)}{n}$ holds for any $m<n$. In paper~[1] it has been proven that the least (exceeding 5040) counterexample to inequality~(\ref{1}) (if exists) is a SA number. Thus, if we can prove the inequality~(\ref{1}) for SA numbers greater than 5040, it will be proven in the general case.

In our numerical experiment, we implement an orderly search of SA numbers, introducing the function $\Omega\colon \mathbb{N} \to \mathbb{N}$, defined by the equality $\Omega(n)=\sum_{i=1}^{k}a_i$. For each value $\omega$ of this function (from 9 to 90), we find the maximum value of the Gr\"onwall function
$$
G_{\max}(\omega)=\max{\{G(n)\mid\Omega(n)=\omega\}}
$$
and mark the number $n_{\omega}^{*}\in\Omega^{-1}(\omega)$, for which this maximum is reached: $G(n_{\omega}^{*})=G_{\max}(\omega)$.

Values $\omega<9$ are of no interest, since the greatest value $\Omega(n)$ for SA numbers~$n$, not exceeding 5040, is equal $\Omega(5040)=\Omega(2^4\cdot 3^2\cdot 5\cdot 7)=8$. At $\omega>90$, the search for number is too cumbersome for us.

The computation results are shown in Table 1. 
In the second column the factorization of $n_{\omega}^{*}$ is written. Dots mean missed primes, each of which occurs in the first power.
The third column shows the number of $n_{\omega}^{*}$ in the sequence of SA-numbers~\cite{Noe}. Each value of $n_{\omega}^{*}$ is followed by its natural logarithm and the least prime number $p(\omega)$ not included in the factorization of $n_{\omega}^{*}$. Ticks in the last column of the table indicate that the following inequality is true
\begin{equation}
\label{4}
p(\omega)<\log n_{\omega}^{*},
\end{equation}
resulting from inequality~(\ref{2}) at $n=n_{\omega}^{*}$ and $p_l=p(\omega)$; exponent $a_l$ is equal to zero. Our calculations show that for the considered values of $\omega$ (at $n=n_{\omega}^{*}$) the inequality~(\ref{2}) holds only if $a_l=0$.

\newpage
\begin{flushright}
Table 1 $\hspace{2.1cm}$
\end{flushright}

\small
\begin{center}

\begin{tabular}{|c|c|c|c|c|c|}

\hline
\par & & & & & \\$\omega$ & $n_{\omega}^{*}$ & SA & $\log n_{\omega}^{*}$ & $p(\omega)$ & Inequality (4)\\
\par & & & & & \\\hline
\par & & & & & \\
9 & $2^5\cdot 3^2 \cdot 5\cdot 7$ & 20 & 9.2 & 11 & \\
10 & $2^5\cdot 3^2 \cdot 5\cdot 7\cdot 11$ & 25 & 11.6 & 13 & \\
11 & $2^5\cdot 3^2 \cdot 5\cdot 7\cdot 11\cdot 13$ & 32 & 14.2 & 17 & \\
12 & $2^5\cdot 3^3 \cdot 5\cdot 7\cdot 11\cdot 13$ & 35 & 15.3 & 17 & \\
13 & $2^5\cdot 3^3 \cdot 5^2 \cdot 7\cdot 11\cdot 13$ & 39 & 16.9 & 17 & \\
14 & $2^5\cdot 3^3 \cdot 5^2 \cdot 7\cdots 17$ & 46 & 19.7 & 19 &\checkmark \\
15 & $2^5\cdot 3^3 \cdot 5^2 \cdot 7\cdots 19$ & 55 & 22.7 & 23 & \\
16 & $2^5\cdot 3^3 \cdot 5^2 \cdot 7\cdots 23$ & 62 & 25.8 & 29 & \\
17 & $2^6\cdot 3^3 \cdot 5^2 \cdot 7\cdots 23$ & 63 & 26.5 & 29 & \\
18 & $2^6\cdot 3^3 \cdot 5^2 \cdot 7\cdots 29$ & 74 & 29.9 & 31 & \\
19 & $2^6\cdot 3^3 \cdot 5^2 \cdot 7\cdots 31$ & 85 & 33.3 & 37 & \\
20 & $2^6\cdot 3^3 \cdot 5^2 \cdot 7^2\cdot 11\cdots 31$ & 91 & 35.2 & 37 & \\
21 & $2^6\cdot 3^4 \cdot 5^2 \cdot 7^2\cdot 11\cdots 31$ & 94 & 36.3 & 37 & \\
22 & $2^6\cdot 3^4 \cdot 5^2 \cdot 7^2\cdot 11\cdots 37$ & 106 & 40.0 & 41 & \\
23 & $2^6\cdot 3^4 \cdot 5^2 \cdot 7^2\cdot 11\cdots 41$ & 116 & 43.7 & 43 &\checkmark \\
24 & $2^6\cdot 3^4 \cdot 5^2 \cdot 7^2\cdot 11\cdots 43$ & 127 & 47.4 & 47 &\checkmark \\
25 & $2^6\cdot 3^4 \cdot 5^2 \cdot 7^2\cdot 11\cdots 47$ & 137 & 51.2 & 53 & \\
26 & $2^7\cdot 3^4 \cdot 5^2 \cdot 7^2\cdot 11\cdots 47$ & 138 & 52.0 & 53 & \\
27 & $2^7\cdot 3^4 \cdot 5^2 \cdot 7^2\cdot 11\cdots 53$ & 149 & 55.9 & 59 & \\
28 & $2^7\cdot 3^4 \cdot 5^2 \cdot 7^2\cdot 11\cdots 59$ & 162 & 60.0 & 61 & \\
29 & $2^7\cdot 3^4 \cdot 5^2 \cdot 7^2\cdot 11\cdots 61$ & 176 & 64.1 & 67 & \\
30 & $2^7\cdot 3^4 \cdot 5^3 \cdot 7^2\cdot 11\cdots 61$ & 181 & 65.7 & 67 & \\
31 & $2^7\cdot 3^4 \cdot 5^3 \cdot 7^2\cdot 11\cdots 67$ & 196 & 69.9 & 71 & \\
32 & $2^7\cdot 3^4 \cdot 5^3 \cdot 7^2\cdot 11\cdots 71$ & 212 & 74.2 & 73 &\checkmark \\
33 & $2^7\cdot 3^4 \cdot 5^3 \cdot 7^2\cdot 11\cdots 73$ & 224 & 78.5 & 79 & \\
34 & $2^7\cdot 3^4 \cdot 5^3 \cdot 7^2\cdot 11^2\cdot 13\cdots 73$ & 231 & 80.9 & 79 &\checkmark \\
35 & $2^7\cdot 3^4 \cdot 5^3 \cdot 7^2\cdot 11^2\cdot 13\cdots 79$ & 246 & 85.3 & 83 &\checkmark \\
36 & $2^7\cdot 3^4 \cdot 5^3 \cdot 7^2\cdot 11^2 \cdot 13\cdots 83$ & 259 & 89.7 & 89 &\checkmark \\
37 & $2^7\cdot 3^4 \cdot 5^3 \cdot 7^2\cdot 11^2\cdot 13\cdots 89$ & 272 & 94.2 & 97 & \\
38 & $2^8\cdot 3^4 \cdot 5^3 \cdot 7^2\cdot 11^2\cdot 13\cdots 89$ & 273 & 94.9 & 97 & \\
39 & $2^8\cdot 3^5 \cdot 5^3 \cdot 7^2\cdot 11^2\cdot 13\cdots 89$ & 276 & 96.0 & 97 & \\

40 & $2^8\cdot 3^5 \cdot 5^3 \cdot 7^2\cdot 11^2\cdot 13\cdots 97$ & 288 & 100.5 & 101 & \\

41 & $2^8\cdot 3^5 \cdot 5^3 \cdot 7^2\cdot 11^2\cdot 13\cdots 101$ & 299 & 105.2 & 103 &\checkmark \\

42 & $2^8\cdot 3^5 \cdot 5^3 \cdot 7^2\cdot 11^2\cdot 13\cdots 103$ & 311 & 109.8 & 107 &\checkmark \\

43 & $2^7\cdot 3^4 \cdot 5^2 \cdot 7^2\cdot 11^2\cdot 13^2\cdot 17\cdots 103$ & 317 & 112.4 & 107 &\checkmark \\

44 & $2^7\cdot 3^5 \cdot 5^3 \cdot 7^2\cdot 11^2\cdot 13^2\cdot 17\cdots 109$ & 341 & 121.0 & 113 &\checkmark \\

45 & $2^7\cdot 3^5 \cdot 5^3 \cdot 7^2\cdot 11^2\cdot 13^2\cdot 17\cdots 113$ & 354 & 125.8 & 127 & \\

46 & $2^8\cdot 3^5 \cdot 5^3 \cdot 7^2\cdot 11^2\cdot 13^2\cdot 17\cdots 113$ & 356 & 126.4 & 127 & \\

47 & $2^8\cdot 3^5 \cdot 5^3 \cdot 7^2\cdot 11^2\cdot 13^2\cdot 17\cdots 127$ & 368 & 131.3 & 131 &\checkmark \\

48 & $2^8\cdot 3^5 \cdot 5^3 \cdot 7^2\cdot 11^2\cdot 13^2\cdot 17\cdots 131$ & 380 & 136.2 & 137 & \\

49 & $2^8\cdot 3^5 \cdot 5^3 \cdot 7^2\cdot 11^2\cdot 13^2\cdot 17\cdots 137$ & 394 & 141.1 & 139 &\checkmark \\

50 & $2^8\cdot 3^5 \cdot 5^3 \cdot 7^2\cdot 11^2\cdot 13^2\cdot 17\cdots 139$ & 408 & 146.0 & 149 & \\

\par & & & & & \\
\hline
\end{tabular}
\end{center}

\large
\newpage

\begin{flushright}
Continue of Table 1 $\hspace{1cm}$
\end{flushright}

\small
\begin{center}
\begin{tabular}{|c|c|c|c|c|c|}

\hline
\par & & & & & \\$\omega$ & $n_{\omega}^{*}$ & SA & $\log n_{\omega}^{*}$ & $p(\omega)$ & Inequality (4)\\
\par & & & & & \\\hline
\par & & & & & \\

51 & $2^9\cdot 3^5 \cdot 5^3 \cdot 7^2\cdot 11^2\cdot 13^2\cdot 17\cdots 139$ & 409 & 146.7 & 149 & \\

52 & $2^8\cdot 3^5 \cdot 5^3 \cdot 7^2\cdot 11^2\cdot 13^2\cdot 17\cdots 151$ & 438 & 156.0 & 157 & \\

53 & $2^9\cdot 3^5 \cdot 5^3 \cdot 7^2\cdot 11^2\cdot 13^2\cdot 17\cdots 151$ & 440 & 156.7 & 157 & \\

54 & $2^8\cdot 3^5 \cdot 5^3 \cdot 7^3\cdot 11^2\cdot 13^2\cdot 17\cdots 157$ & 458 & 163.0 & 163 &\checkmark \\

55 & $2^9\cdot 3^5 \cdot 5^3 \cdot 7^3 \cdot 11^2\cdot 13^2\cdot 17\cdots 157$ & 459 & 163.7 & 163 &\checkmark \\

56 & $2^9\cdot 3^5 \cdot 5^3 \cdot 7^3\cdot 11^2\cdot 13^2\cdot 17\cdots 163$ & 476 & 168.9 & 167 &\checkmark \\

57 & $2^9\cdot 3^5 \cdot 5^3 \cdot 7^3\cdot 11^2\cdot 13^2\cdot 17\cdots 167$ & 493 & 173.9 & 173 &\checkmark \\

58 & $2^8\cdot 3^5 \cdot 5^3 \cdot 7^3\cdot 11^2\cdot 13^2\cdot 17^2\cdot 19\cdots 173$ & 518 & 181.2 & 179 &\checkmark \\

59 & $2^9\cdot 3^5 \cdot 5^3 \cdot 7^3\cdot 11^2\cdot 13^2\cdot 17^2\cdot 19\cdots 173$ & 519 & 181.9 & 179 &\checkmark \\

60 & $2^8\cdot 3^5 \cdot 5^3 \cdot 7^3\cdot 11^2\cdot 13^2\cdot 17^2\cdot 19\cdots 181$ & 554 & 191.6 & 191 &\checkmark \\

61 & $2^9\cdot 3^5 \cdot 5^3 \cdot 7^3\cdot 11^2\cdot 13^2\cdot 17^2\cdot 19\cdots 181$ & 555 & 192.3 & 191 &\checkmark \\

62 & $2^9\cdot 3^5 \cdot 5^3 \cdot 7^3\cdot 11^2\cdot 13^2\cdot 17^2\cdot 19\cdots 191$ & 575 & 197.8 & 193 &\checkmark \\

63 & $2^9\cdot 3^5 \cdot 5^3 \cdot 7^3\cdot 11^2\cdot 13^2\cdot 17^2\cdot 19\cdots 193$ & 596 & 202.8 & 197 &\checkmark \\

64 & $2^9\cdot 3^5 \cdot 5^3 \cdot 7^3\cdot 11^2\cdot 13^2\cdot 17^2\cdot 19\cdots 197$ & 613 & 208.1 & 199 &\checkmark \\

65 & $2^9\cdot 3^5 \cdot 5^3 \cdot 7^3\cdot 11^2\cdot 13^2\cdot 17^2\cdot 19\cdots 199$ & 628 & 213.4 & 211 &\checkmark \\

66 & $2^9\cdot 3^5 \cdot 5^3 \cdot 7^3\cdot 11^2\cdot 13^2\cdot 17^2\cdot 19\cdots 211$ & 643 & 218.8 & 223 & \\

67 & $2^9\cdot 3^5 \cdot 5^3 \cdot 7^3\cdot 11^2\cdot 13^2\cdot 17^2\cdot 19^2\cdot 23\cdots 211$ & 653 & 221.7 & 223 & \\

68 & $2^9\cdot 3^5 \cdot 5^3 \cdot 7^3\cdot 11^2\cdot 13^2\cdot 17^2\cdot 19^2\cdot 23\cdots 223$ & 670 & 227.1 & 227 &\checkmark \\

69 & $2^9\cdot 3^5 \cdot 5^3 \cdot 7^3\cdot 11^2\cdot 13^2\cdot 17^2\cdot 19^2\cdot 23\cdots 227$ & 685 & 232.5 & 229 &\checkmark \\

70 & $2^9\cdot 3^5 \cdot 5^3 \cdot 7^3\cdot 11^2\cdot 13^2\cdot 17^2\cdot 19^2\cdot 23\cdots 229$ & 701 & 238.0 & 233 &\checkmark \\

71 & $2^9\cdot 3^5 \cdot 5^3 \cdot 7^3\cdot 11^2\cdot 13^2\cdot 17^2\cdot 19^2\cdot 23\cdots 233$ & 717 & 243.4 & 239 &\checkmark \\

72 & $2^9\cdot 3^5 \cdot 5^3 \cdot 7^3\cdot 11^2\cdot 13^2\cdot 17^2\cdot 19^2\cdot 23\cdots 239$ & 733 & 248.9 & 241 &\checkmark \\

73 & $2^9\cdot 3^5 \cdot 5^3 \cdot 7^3\cdot 11^2\cdot 13^2\cdot 17^2\cdot 19^2\cdot 23\cdots 241$ & 748 & 254.4 & 251 &\checkmark \\

74 & $2^9\cdot 3^5 \cdot 5^4 \cdot 7^3\cdot 11^2\cdot 13^2\cdot 17^2\cdot 19^2\cdot 23\cdots 241$ & 752 & 256.0 & 251 &\checkmark \\

75 & $2^9\cdot 3^6 \cdot 5^4 \cdot 7^3\cdot 11^2\cdot 13^2\cdot 17^2\cdot 19^2\cdot 23\cdots 241$ & 755 & 257.1 & 251 &\checkmark \\

76 & $2^9\cdot 3^6 \cdot 5^4 \cdot 7^3\cdot 11^2\cdot 13^2\cdot 17^2\cdot 19^2\cdot 23\cdots 251$ & 774 & 262.6 & 257 &\checkmark \\

77 & $2^9\cdot 3^6 \cdot 5^4 \cdot 7^3\cdot 11^2\cdot 13^2\cdot 17^2\cdot 19^2\cdot 23\cdots 257$ & 791 & 268.2 & 263 &\checkmark \\

78 & $2^9\cdot 3^6 \cdot 5^4 \cdot 7^3\cdot 11^2\cdot 13^2\cdot 17^2\cdot 19^2\cdot 23\cdots 263$ & 808 & 273.7 & 269 &\checkmark \\

79 & $2^9\cdot 3^6 \cdot 5^4 \cdot 7^3\cdot 11^2\cdot 13^2\cdot 17^2\cdot 19^2\cdot 23\cdots 269$ & 825 & 279.3 & 271 &\checkmark \\

80 & $2^9\cdot 3^6 \cdot 5^4 \cdot 7^3\cdot 11^2\cdot 13^2\cdot 17^2\cdot 19^2\cdot 23\cdots 271$ & 842 & 284.9 & 277 &\checkmark \\

81 & $2^9\cdot 3^6 \cdot 5^4 \cdot 7^3\cdot 11^2\cdot 13^2\cdot 17^2\cdot 19^2\cdot 23\cdots 277$ & 859 & 290.6 & 281 &\checkmark \\

82 & $2^9\cdot 3^6 \cdot 5^4 \cdot 7^3\cdot 11^2\cdot 13^2\cdot 17^2\cdot 19^2\cdot 23\cdots 281$ & 874 & 296.2 & 283 &\checkmark \\

83 & $2^9\cdot 3^6 \cdot 5^4 \cdot 7^3\cdot 11^2\cdot 13^2\cdot 17^2\cdot 19^2\cdot 23\cdots 283$ & 889 & 301.8 & 293 &\checkmark \\

84 & $2^9\cdot 3^6 \cdot 5^4 \cdot 7^3\cdot 11^2\cdot 13^2\cdot 17^2\cdot 19^2\cdot 23\cdots 293$ & 903 & 307.5 & 307 &\checkmark \\

85 & $2^{10}\cdot 3^6 \cdot 5^4 \cdot 7^3\cdot 11^2\cdot 13^2\cdot 17^2\cdot 19^2\cdot 23\cdots 293$ & 904 & 308.2 & 307 &\checkmark \\

86 & $2^{10}\cdot 3^6 \cdot 5^4 \cdot 7^3\cdot 11^2\cdot 13^2\cdot 17^2\cdot 19^2\cdot 23^2 \cdot29 \cdots 293$ & 912 & 311.4 & 307 &\checkmark \\

87 & $2^{10}\cdot 3^6 \cdot 5^4 \cdot 7^3\cdot 11^2\cdot 13^2\cdot 17^2\cdot 19^2\cdot 23^2 \cdot29 \cdots 307$ & 927 & 317.1 & 311 &\checkmark \\

88 & $2^{10}\cdot 3^6 \cdot 5^4 \cdot 7^3\cdot 11^2\cdot 13^2\cdot 17^2\cdot 19^2\cdot 23^2 \cdot29 \cdots 311$ & 942 & 322.8 & 313 &\checkmark \\

89 & $2^9\cdot 3^6 \cdot 5^4 \cdot 7^3\cdot 11^2\cdot 13^2\cdot 17^2\cdot 19^2\cdot 23^2 \cdot29 \cdots 317$ & 971 & 333.6 & 331 &\checkmark \\

90 & $2^{10}\cdot 3^6 \cdot 5^4 \cdot 7^3\cdot 11^2\cdot 13^2\cdot 17^2\cdot 19^2\cdot 23^2 \cdot29 \cdots 317$ & 972 & 334.3 & 331 &\checkmark \\

\par & & & & & \\
\hline
\end{tabular}
\end{center}
\large
\newpage

\section{Concluding remarks}

From Table 1 we can see that at the increase in parameter $\omega$, the proportion of its values $\omega'\le\omega$, for which the inequality~(\ref{4}) is satisfied, rises (which, however, does not rule out a fundamental change in the situation at greater values of $\omega$). If we can prove that
$$
\exists\tilde\omega\qquad\mbox{s. t.}\qquad \forall\omega\in[9,\tilde\omega]\hspace{0.5cm}G_{\max}<e^\gamma\qquad
\mbox{and}\qquad
\forall\omega\ge\tilde\omega \hspace{0.5cm}p(\omega)<\log n_{\omega}^{*},
$$
we get a proof of the Riemann hypothesis.

In fact, the fulfillment of the inequality~(\ref{4}) for certain $\omega$ value means that the inequality $G_{\max}(\omega+1)>G_{\max}(\omega)$ also holds:
$$
G_{\max}(\omega+1)=G(n_{\omega+1}^{*})\ge
G(n_{\omega}^{*}p(\omega))>G_{\max}(\omega).
$$
Let us consider any SA number $n$ such that $\Omega(n)\ge\tilde\omega$. If we construct an infinite sequence
$$
u_1=n,\mbox{ }u_2=n_{\Omega(n)+1}^{*},
\mbox{ }u_3=n_{\Omega(n)+2}^{*},\mbox{ ... , }
u_i=n_{\Omega(n)+i-1}^{*},\mbox{ ...}
$$
we will note that it has the following property:
$$
\forall i,j\qquad i<j \Longrightarrow G(u_i)<G(u_j).
$$
Given the Gr\"onwall equality~\cite{Gronwall}
$$
\limsup\limits_{n \rightarrow \infty} G(n)=e^\gamma,
$$
we can state that Robin's inequality~(\ref{1}) is true for all SA numbers that exceed 5040, and thus holds in the general case, which is equivalent to the Riemann hypothesis.

\end{document}